# A Robust Optimization Approach to a Real Humanitarian Cold Supply Chain Planning on the COVID-19 crisis

**Behnam Malmir and Christopher W. Zobel**

**Department of Business Information Technology,**

**Pamplin College of Business, Virginia Tech, Blacksburg, Virginia, USA**

**Abstract**

In this study, a vaccine supply chain model is developed considering the humanitarian aspects under uncertain conditions in Iran. There are three main components to supply the required vaccines for vaccination centers that can be designed and managed within three echelons: suppliers, distribution centers, and vaccination centers. Iran's Ministry of Health and Medical (IMHM), Iranian Red Crescent Society (IRCS), and private sector companies play the role of suppliers in the model. A robust optimization approach is employed to address real-world uncertainty and find a solution dealing with all uncertain data possibilities and equity of the supply chain. To this aim, the actual COVID-19 data of Iran is gathered, including data on five different types of vaccines used in Iran. Then, it is investigated how vaccination programs can be accomplished more efficiently by considering priority issues, transportation, holding costs, and deprivation costs in Tehran, Iran's largest city and one of the most populous cities globally. Next, a set of sensitivity analyses are conducted, comprising assessing the variation of deprivation costs for different problem instances as well as the budgets assigned to each supplier. One of the critical assumptions in this research is prioritizing older-aged groups over younger-aged ones into account for receiving the vaccine.

**Keywords:** Humanitarian supply chain; Vaccination; Robust Optimization, Uncertainty; Equity; COVID-19.

## 3. Problem description

In this study, a vaccine supply chain model is proposed by considering the humanitarian aspects. In this regard, there are three main components to supply the vaccines for vaccination centers that can be designed and managed in three echelons: suppliers, distribution centers, and vaccination centers.

### 3.1 Assumptions

The main assumptions considered in modeling the on-hand humanitarian supply chain problem can be described as follows:

i. Iran's Ministry of Health and Medical (IMHM), Iranian Red Crescent Society (IRCS), and private sector companies play the role of suppliers in the considered model. In other words, any other source of vaccine supply must be directed in one of these channels.

ii. Suppliers can provide more than one type of vaccine.

*iii.* There is a set of predefined locations for the distribution centers and vaccination centers.

*iv.* Each distribution center can serve more than one vaccination center.

*v.* The suppliers have capacity limitations for providing each type of vaccine, while all of them can play this role.

*vi.* Vaccine distribution flow is allowed only from suppliers to distribution centers and from them to the vaccination centers; it is not allowed to directly distribute from suppliers to vaccination centers. In other words, it is assumed that direct vaccine distribution flows cannot exist from one supplier inventory to another or from one distribution center to another.

*vii.* No limitation is assumed for the number of vehicles. Furthermore, vaccination centers are accessible through the current road network.

*viii.* Ten age groups are taken into account including …

*ix.* A minimum demand satisfaction rate is defined for each type of vaccine in any vaccination center, giving the minimum service level in vaccination centers with the goal of satisfying equity.

*x.* The daily demand for vaccines in each vaccination center is uncertain.

*xi.* The planning horizon for the problem to get interventions is known.

*xii.* Those parts of vaccines which cannot be satisfied will be backordered.

The following depiction provides an overview of the health system entities that contributed to the supply of vaccines. As mentioned, the proposed model makes some additional assumptions about the vaccine supply chain. For instance, although there can be a direct connection between supplier inventory and vaccination centers, it is assumed that vaccines are generally transported from distribution centers to vaccination centers. In an emergency condition, the government can provide essential resources such as financial, transportation, and inventory infrastructure. Hence, future consideration of improving suppliers' inventory might enhance the effectiveness of delivering vaccines in proper time, leading to a reduction in backorder demands.

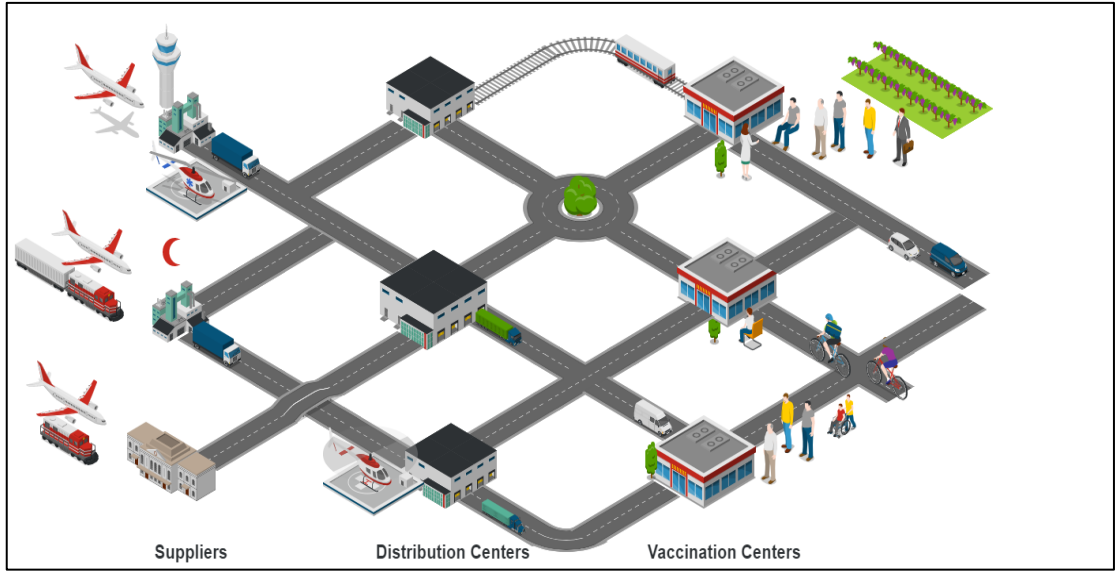

**Fig. 1.** Logistic network for the considered three-echelon vaccine supply chain.

### 3.2 Notations

#### 3.2.1 Sets and indices

$i$: Index for suppliers ($i =1, 2, 3$),

$j$: Index for distribution center ($j =1, 2,…, J$),

$k$: Index for vaccination center ($k =1, 2,..., K$),

$v$: Index for vaccine types ($v =1, 2,..., V$),

$a$: Index for age groups ($a =1, 2,…, A$),

$t$: Index for time periods ($t =1, 2,…,T$).

#### 3.2.2 Input parameters

$l_i$: Lead time for supplier $i$,

$price_{iv}$: Price of vaccine type $v$ purchased from supplier $i$,

$Bud_i$: Available budget to purchase vaccines from supplier $i$,

$Mxord_{iv}$: Maximum order of vaccine type $v$ placed to supplier $i$,

$Cap_i$: Maximum capacity of supplier $i$,

$Capd_j$: Maximum capacity of distribution center $j$,

$CostI_{ijv}$: Transportation cost for vaccine type $v$ from supplier $i$ to distribution center $j$,

$CI_{ij}$: Transportation fixed cost from supplier $i$ to distribution center $j$,

$Costd_{jkv}$: Transportation cost for vaccine type $v$ from distribution center $j$ to vaccine center $k$,

$Cd_{jk}$: Transportation fixed cost from distribution center $j$ to vaccination center $k$,

$h_v$: Holding cost of vaccine type $v$ for one period,

$\theta_f$: Weight of objective function $f$,

$r(t)$: Deprivation intensity function at period $t$,

$\omega$: Minimum percentage of demand for vaccine $v$ in vaccination center $k$ that should be satisfied in period $t$,

$inv_{ivt}$: Inventory level of vaccine type $v$ in supplier $i$ at period $t$,

$invd_{jvt}$: Inventory level of vaccine type $v$ in distribution center $j$ at period $t$,

$t_v$: Maximum allowed holding time for vaccine type $v$.

$d_{kvt}$: Demand for vaccine $v$ in vaccination center $k$ at period $t$,

$Demand_{aj}$: Total demand of age group $a$ in distribution center $j$,

### 3.2.3 Decision variables

$X_{ivt}$: Amount of ordered vaccine type $v$ to supplier $i$ at period $t$,

$Y_{ijvt}$: Amount of vaccine type $v$ transported from supplier $i$ to distribution center $j$ at period $t$,

$W_{jkvt}$: Amount of vaccine type $v$ transported from distribution center $j$ to vaccination center $k$.

$u_{at}$: Priority of age group $a$ at period $t$; a binary variable,

$u_{at}$: Equals 1 if age group $a$ receive vaccine at period $t$; a binary variable,

$back_{kvt}$: Amount of back order for vaccine type $v$ in vaccination center $k$ at period $t$,

$q_{kvt}$: Amount of demand for vaccine type $v$ in vaccination center $k$ that is satisfied at period $t$,

$yb_{ijt}$: Equals 1 if there are items to transport from supplier $i$ to distribution center $j$ at period $t$, otherwise 0,

$wb_{jkt}$: Equals 1 if there are items to transport from distribution center $j$ to vaccination center $k$ at period $t$, otherwise 0,

The proposed nonlinear mathematical model is formulated as follows:

### 3.3 The mathematical model

| | |
|---|---|
| minimize $Z = \theta_1 P_1 + \theta_2 P_2 + \theta_3 P_3$ | (1) |
| $P_1 = \sum_i \sum_v \sum_t inv_{ivt} \times h_v + \sum_j \sum_v \sum_t invd_{jvt} \times h_v$ | (2) |
| $P_2 = \sum_i \sum_j \sum_v \sum_t (CI_{ij} \times yb_{ijt} + Y_{ijvt} \times CostI_{ijv})$ $+ \sum_j \sum_k \sum_v \sum_t (Cd_{jk} \times wb_{jkt} + W_{jkvt} \times Costd_{jkv})$ | (3) |
| $P_3 = \sum_t \sum_k \sum_v r(t) \times (q_{kvt} - q'_{kvt})$ | (4) |
| subject to | |
| $X_{ivt} \leq maxord_{iv} \quad \forall i, t,$ | (5) |
| $\sum_v \sum_t price_{iv} \times X_{ivt} \leq Bud_i \quad \forall i,$ | (6) |
| $\sum_v \sum_t X_{ivt} - \sum_j \sum_v \sum_t Y_{ijvt} \leq Cap_i \quad \forall i,$ | (7) |
| $Y_{ijvt} \leq \sum_{t \in [t-t_v, t-l_i]} X_{ivt} \quad \forall i, j, v, t,$ | (8) |
| $inv_{ivt} = \sum_{t \in [t-t_v, t-l_i]} X_{ivt} - \sum_j \sum_t Y_{ijvt} \quad \forall i, v, t,$ | (9) |
| $\sum_v inv_{ivt} \leq Cap_i \quad \forall i, t,$ | (10) |
| $invd_{jvt} = \sum_i \sum_{t \in [t-t_v, t]} Y_{ijvt} - \sum_k \sum_{t \in [t-t_v, t]} W_{jkvt} \quad \forall j, v, t$ | (11) |
| $\sum_v invd_{jvt} \leq Capd_j \quad \forall j, t$ | (12) |
| $q_{kvt} = d_{kvt} + back_{kvt-1} \quad \forall k, v, t$ | (13) |
| $back_{kvt} = q_{kvt} - \sum_j W_{jkvt} \quad \forall k, v, t$ | (14) |
| $\omega_{kvt} \times q_{kvt} \leq \sum_j W_{jkvt} \quad \forall k, v, t$ | (15) |
| $Q_{jvt} = \sum_{k \in j_k} q_{kvt} \quad \forall j, v, t$ | (16) |
| $Q'_{jvt} = \sum_{k \in K} q'_{kvt} \quad \forall j, v, t$ | (17) |
| $\left\|\dfrac{Q'_{jvt}}{\sum_{a\in[1,10]} Demand_{aj}} - \dfrac{Q'_{j'vt}}{\sum_{a\in[1,10]} Demand_{aj'}}\right\| \leq \xi \qquad \forall j, j', v, t$ | (18) |
| $\sum_{k\in K}\sum_{v} q_{kvt} \leq \sum_{a\in[1,10]} u_{at} \times Demand_{aj} \qquad \forall j, t$ | (19) |
| $u_{a+1.t} \leq u_{at} \qquad \forall i, t,$ | (20) |
| $Y_{ijvt} \leq yb_{ijt} \times M$ | (21) |
| $W_{jkvt} \leq wb_{jkt} \times M$ | (22) |
| $X_{ivt}, Y_{ijvt}, W_{jkvt} \geq 0, u_{at} \in \{0,1\} \qquad \forall i, j, v, k, t.$ | (23) |

The proposed objective function, Eq. (1), contains three terms, including the imposed costs to the system that should be minimized, as follows: the first part ($P_1$) includes total holding costs in the predefined network, the second part ($P_2$) comprises total transportation costs, and the third part ($P_3$) calculates total deprivation costs. As previously mentioned, minimization of deprivation costs is an approach that this article proposes to ensure vertical equity in the under-study society.

It is noteworthy that the summation of $\theta_1, \theta_2$ and $\theta_3$ is 1. The weight or importance of the first, second, and third parts of the objective function are considered as user-defined values $\theta_1, \theta_2$ and $\theta_3$. In order to show the importance of the social component in the humanitarian context, in this study, we specify values as $\theta_1 = 0.3, \theta_2 = 0.1$ and $\theta_3 = 0.6$. In the transportation costs, both variable and fixed costs for vehicles are incorporated by $CI_{ij}$, $Cd_{jk}$, and $CostI_{ijv}$ and $Costd_{jkv}$ respectively.

Supposing that $\lambda_k$ be the time interval at which vaccines are distributed to a vaccination center $k$, one can denote the deprivation cost through the deprivation intensity function $r(t)$. We take the summation of all charges over the deliveries. It should be noted that $r(t)$, defined as $r(t) = 3t$ in this study (Malmir and Zobel, 2021), takes the value of zero when the sufficient supply delivers to a vaccination center. Researchers took varying forms of the $r(t)$ function; for instance, Gutjahr and Fischer (2018) have taken quadratic functions. In this article, a simple linear approximation of the quadratic function was considered, according to Malmir and Zobel (2021). Alternative forms of deprivation function and the related effects on the performance of the proposed system might be the subject of future research studies.

Constraint (5) shows that based on the level of international relationships of the suppliers and capability of importing goods, only the restricted amounts of each type of vaccine can be ordered. Constraint (6) represents the budget limitation for suppliers in their decisions. The constraint associated with the restricted capacity of suppliers is stated in the form of constraint (7). Constraint (8) demonstrates the flow of supplied vaccines between suppliers and distribution centers. Calculation of each supplier inventory level for different types of vaccines and capacity limitation for each inventory is presented in Constraints (9) and (10), respectively. The level inventory for each type of vaccine in distribution centers is calculated in Constraint (11). Constraint (12) shows that any distribution center based on its limited equipment has a limited capacity for the overall inventory of vaccines. Constraints (13) and (14) show that in each period, the quantity of demand is the summation of the amount of backordered demand from previous periods and the new demand arisen in the current period. Constraint (15) is related to the vertical equity and shows that a minimum rate of demand should be satisfied for each vaccination center. The total quantity of demand in each distribution center and the demand to be satisfied for each distribution center are presented in Constraints (16) and (17), respectively. Constraint (18) is the equity constraint that ensures the fair and equitable distribution of resources. Constraint (19) is related to age groups and ensures that the total demands are equal to or less than the targeted population of prioritized age groups. Moreover, Constraint (20) shows that older age groups have priority over younger age ones. If the associated transportations exist, logical constraints (21) and (22) ensure that the binary variables will return 1.

## 4. Robust Optimization

In recent years, researchers in optimization have paid a great deal of attention to uncertain data in optimization problems. In this regard, varying techniques are proposed, such as stochastic programming (SP), robust optimization (RO), and Distributionally Robust Optimization (DRO). Robust optimization attempts to find a solution dealing with all uncertain data possibilities. Researchers have proposed different robust optimization methodologies. Soyester (1973) proposed the first robust model for linear optimization models, including uncertain data. Bertsimas and Sim (2004) developed a robust formulation investigating parametric data uncertainty without overly penalizing the objective function. Bertsimas and Sim (2004) controlled the constraint violations probability and the consequence of conservatism on the objective function. This article develops a Bertsimas' and Sim's optimization method to solve the proposed vaccine supply chain model.

Considering the following mathematical problem, which includes a set of $j$ variables in Eq. (2):

minimize $C^T x$ (22)

subject to (23)

$Ax \leq b,$ (24)

$lb \leq x \leq ub.$ (25)

Uncertainty in some parameters affects both matrix $\boldsymbol{A} = (a_{n,j})$ and vector $\boldsymbol{C}$. These features are modeled as symmetric and restricted random variables such as $\tilde{c}_n$ and $\tilde{a}_n, j \in J_n$. Thus, another form for defining uncertain parameters through "mean value" and "range of each uncertain factor" are is follows:

$\tilde{c}_n \in [c_n - \hat{c}_n, c_n + \hat{c}_n] \quad \forall \tilde{c}_n \in C,$ (26)

$\tilde{a}_{nj} \in [a_{nj} - \hat{a}_{nj}, a_{nj} + \hat{a}_{nj}] \quad \forall \tilde{a}_{nj} \in A.$ (27)

Here, $\Gamma_n$ ($n = 0,1, \ldots, CN$) as conversion level (CL) or budget of uncertainty is presented for robustness purposes and regulating the model robustness giving the solution that can take different values in the interval of $[0, |J_n|]$, with $J_n$ as a set comprising uncertain factor of the constraint $n$ ($n = 1, \ldots, CN$), $J_n = \{j | \hat{a}_{nj} > 0\}$, or the objective function ($n = 0$). Thus, we can present the value of conversion level for the objective function and constraining $n$ by $\Gamma_0$ and $\Gamma_n$, respectively. Finally, the robust equivalent model could be presented as follows:

$$\text{minimize } c^T x + \max_{\{S_0 \cup \{t_0\} S_0 \subseteq J_0, |S_0| = \lfloor \Gamma_0 \rfloor, t_0 \in J_0 \setminus S_0\}} \left\{ \sum_{j \in S_0} \hat{c}_n |x_j| + (\Gamma_0 - \lfloor \Gamma_0 \rfloor) \hat{c}_{t_0} |x_0| \right\} \quad (28)$$

subject to

$$\sum_j a_{n,j} x_j + \max_{\{S_0 \cup \{t_0\} S_0 \subseteq J_0, |S_0| = \lfloor \Gamma_0 \rfloor, t_0 \in J_0 \setminus S_0\}} \left\{ \sum_{j \in S_0} \hat{a}_{n,j} |x_j| + (\Gamma_n - \lfloor \Gamma_n \rfloor) \hat{a}_{n,t_n} |x_{t_n}| \right\} \leq b_n, \quad (29)$$

$$lb_j \leq x_j \leq ub_j \quad (j = 1,2, \ldots, DN). \quad (30)$$

With the use of duality theory and applying to the above equations, an equivalent optimization model is obtained as follows:

$$\text{minimize } \sum_j c_j x_j + z_0 \Gamma_0 + \sum_{j \in J_0} p_{0j} \quad (31)$$

$$\sum_j a_{nj} x_j + z_n \Gamma_n + \sum_{j \in J_n} p_{nj} \leq b_n \quad (n = 1, \ldots, CN), \quad (32)$$

$$z_0 + p_{0j} \geq \tilde{c}_j y_j \quad \forall j \in J_0, \quad (33)$$

$$z_n + p_{nj} \geq \tilde{a}_{nj} y_j \quad \forall j \in J_n; (n = 1, \ldots, CN), \quad (34)$$

$$-y_j \leq x_j \leq y_j \quad (j = 1, \ldots, CN), \quad (35)$$

$$lb_j \leq x_j \leq ub_j \quad (j = 1, \ldots, CN), \quad (36)$$

$$z_n \geq 0 \quad (j = 1, \ldots, CN), \quad (37)$$

$$P_{nj} \geq 0 \quad \forall j \in J_n; (n = 1, \ldots, CN), \quad (38)$$

$$y_j \geq 0 \quad (j = 1, \ldots, CN). \quad (39)$$

Now, based on the explanations mentioned above, the robust model of the proposed mathematical model in the vaccine supply chain could be expanded by taking a set of symmetric restricted intervals into account to show the uncertainty in the maximum amount of ordering for suppliers and the available budget for each supplier. Accordingly, the final form of the LP model can be presented as follows:

$$\text{minimize } Z = \theta_1 P_1 + \theta_2 P_2 + \theta_3 P_3 \quad (40)$$

$$P_1 = \sum_i \sum_v \sum_t inv_{ivt} \times h_v + \sum_j \sum_v \sum_t invd_{jvt} \times h_v \tag{41}$$

$$P_2 = \sum_i \sum_j \sum_v \sum_t (CI_{ij} + Y_{ijvt} \times CostI_{ijv}) + \sum_j \sum_k \sum_v \sum_t (Cd_{jk} + W_{jkvt} \times Costd_{jkv}) \tag{42}$$

$$P_3 = \sum_t \sum_k \sum_v r(t) \times (q_{kvt} - q'_{kvt}) \tag{43}$$

subject to

$$X_{ivt} + r1_{it}\Gamma + \varphi(H1_{it}) \leq maxord_{iv} \quad \forall i,t, \tag{44}$$

$$\sum_v \sum_t price_{iv} \times X_{ivt} + r2_i\Gamma + \varphi(H2_i) \leq Bud_i \quad \forall i, \tag{45}$$

$$r1_{i,t} + H1_{it} \geq \widehat{maxord}_{iv} \quad \forall i,v, \tag{46}$$

$$r2_i + H2_i \geq \widehat{Bud}_i \quad \forall i, \tag{47}$$

$$r2_i, H2_i, r1_{it}, H1_{it} \geq 0 \quad \forall i,t, \tag{48}$$

Eqs. (7)-(21).

## 5. Experimental Results

### 5.1. Case study

According to official statistics, approximately 5.8 million of Iranian people have been infected with COVID-19, of which nearly 125K individuals have died (Worldometer, 2021). Among all 31 provinces in Iran, Tehran is the most populated province, with almost 13.2 million. Thus, the local government agencies are more concerned about managing the outbreak through managerial decisions on social distancing, lockdowns, and vaccination programs. Although nearly 23% of Iranian people have been fully vaccinated till October 2021 (Google News, 2021), the city of Tehran has been almost labeled as a 'red city' for an extended period. In this case study, especially in Tehran, we are focused on investigating how vaccination programs can be accomplished more efficiently by considering priority issues, transportation and holding costs, as well as deprivation costs.

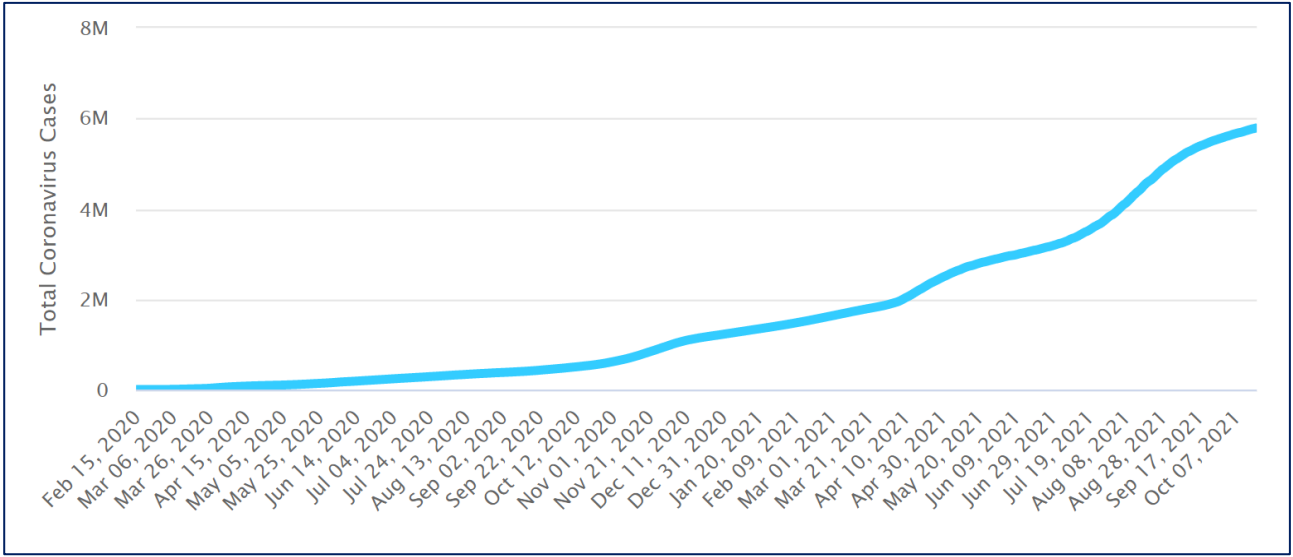

**Fig. 2.** Number of total Coronavirus Cases in Iran from 2020 to 2021.

Moreover, Iran began vaccination in February 2021 with a limited number of imported Sputnik vaccines (WSJ, 2021). Although Iran has started to import vaccines from different countries, five types of vaccines, i.e., Barkat, Sputnik, Novavax, PastoCoVac, and Fakhravak, have been produced so far by Iranian scientists (EuroNews, 22 Oct. 2021). In this study, taking different situations into account, including brands, manufacturer country of vaccine, lead time, and channel through which they were imported, altogether 15 kinds of various vaccines are considered.

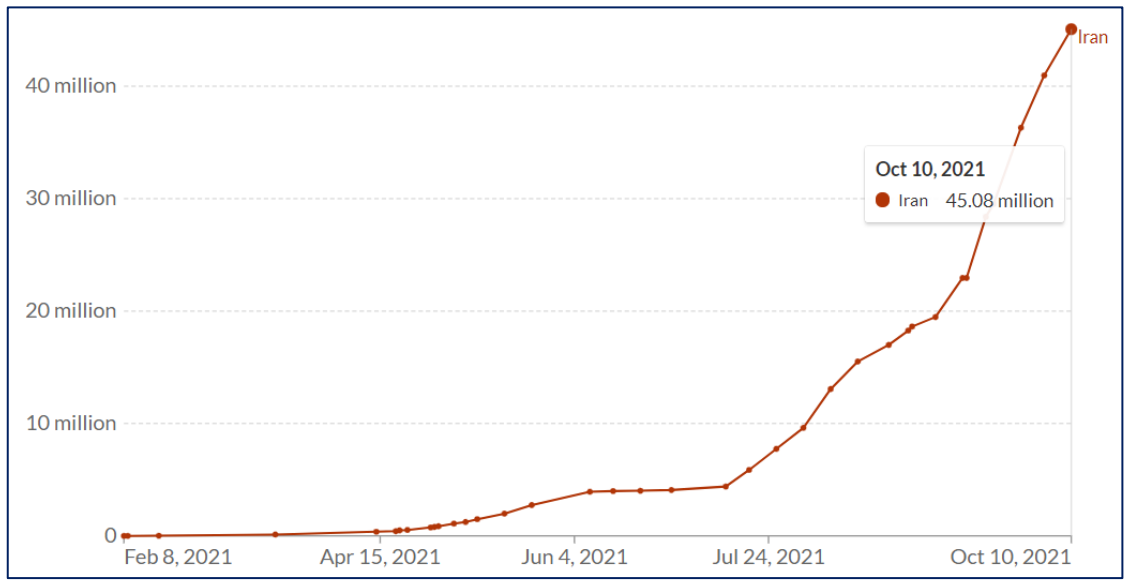

**Fig. 3.** Number of Iranian people who received at least one dose of COVID-19 vaccine.

**Table 1.** Information of conducted experiments.

| #Problem instance | No. of Periods | No. of Distribution Centers | No. of Vaccination Centers | No. of Vaccine Types |
|---|---|---|---|---|
| 1 | 5 | 10 | 20 | 2 |
| 2 | 5 | 10 | 30 | 2 |
| 3 | 10 | 10 | 40 | 2 |
| 4 | 10 | 10 | 40 | 3 |
| 5 | 15 | 10 | 40 | 3 |
| 6 | 20 | 20 | 50 | 4 |
| 7 | 20 | 20 | 50 | 5 |
| 8 | 25 | 20 | 50 | 6 |
| 9 | 25 | 20 | 60 | 8 |
| 10 | 30 | 20 | 60 | 10 |
| 11 | 35 | 31 | 70 | 10 |
| 12 | 40 | 31 | 80 | 10 |
| 13 | 45 | 31 | 85 | 10 |
| 14 | 50 | 31 | 90 | 10 |
| 15 | 100 | 31 | 100 | 15 |

Table 1 represents the information of these problem instances including small (the first five instances), medium (the second five instances), and large size (the last five instances) for experiments.

There are three major vaccine suppliers in Iran, including the IRCS, the private sector, and the IMHM. These suppliers, individually and based on government regulations, have started to import vaccines. Although the suppliers have limitations in terms of available budget, due to existing sanctions against Iran, the capacity of inventories, and the restricted maximum amount of ordering vaccines, they are trying to improve their internal infrastructures and make efficient strategies in the uncertain situation of the COVID-19 outbreak.

Table 2 presents some parameters related to suppliers. It should be noted that since there are assumed different types of vaccines, we have varying maximum ordering for each type of vaccine, and the average maximum order represented in Table 2.

**Table 2.** Parameters related to suppliers.

| Name of Source | Available Budget ($) | Average of Maximum Order (doses) | The capacity of supplier Inventory (doses) |
|---|---|---|---|
| IRCS | 250 Million | 3.18 Million | 3.0 Million |
| Private sector | 400 Million | 3.35 Million | 2.5 Million |
| IMHM | 700 Million | 3.35 Million | 3.5 Million |

**5.2. Results**

To assess and show the performance of the proposed model, different numerical samples are randomly generated and solved using GAMS 23.5 on a Core i5 with 1 GB RAM.

By solving the proposed mathematical model with the mentioned parameters, the objective function value and the computational time are shown in Table 3. Table 3 represents the objective function components $P_1, P_2$, and $P_3$ along with $P^*$ (value of the total objective function). Based on these results, one can initially detect that growing the problem size increases the objective function. Clearly, it is based on the increase in the vaccine demands in different age groups and the rise in the number of vaccination centers. In addition, it could be inferred from the obtained values that the growth rate of deprivation costs decreases while the growth rate of transportation and holding costs increases. This phenomenon can be construed as follows: due to the rise in demands, the proposed model tries to reduce the deprivation costs by distributing more vaccine doses for all vaccination centers and fair distribution. The model presents that after a crisis, such as the COVID-19 pandemic, the minimization of deprivation costs has priority over other expenses like transportation and holding costs.

**Table 3.** Robust model solution.

| #Problem instance | Computational Time (seconds) | $P_1$ | $P_2$ | $P_3$ | Objective Function Value ($) |
|---|---|---|---|---|---|
| 1 | 13 | 223,247 | 446,643 | 232,299 | 902,189 |
| 2 | 18 | 281,243 | 523,199 | 594,996 | 1,399,438 |
| 3 | 25 | 842,316 | 1,579,954 | 1,053,461 | 3,475,731 |
| 4 | 28 | 635,943 | 1,221,070 | 1,640,258 | 3,497,271 |
| 5 | 35 | 1,001,141 | 1,887,835 | 2,293,312 | 5,182,288 |
| 6 | 65 | 1,884,663 | 3,502,917 | 3,310,674 | 8,698,254 |
| 7 | 70 | 1,650,137 | 3,119,714 | 3,941,513 | 8,711,364 |
| 8 | 81 | 2,289,756 | 4,274,120 | 4,221,590 | 10,785,466 |
| 9 | 103 | 2,939,661 | 5,473,503 | 4,601,078 | 13,014,242 |
| 10 | 127 | 3,729,784 | 6,983,337 | 4,864,448 | 15,577,569 |
| 11 | 163 | 5,526,398 | 10,277,000 | 5,430,885 | 21,234,283 |
| 12 | 186 | 7,408,152 | 13,801,007 | 6,431,385 | 27,640,544 |
| 13 | 208 | 8,305,952 | 15,459,416 | 9,336,947 | 33,102,315 |
| 14 | 222 | 10,056,322 | 18,684,245 | 10,141,944 | 38,882,511 |
| 15 | 467 | 24,907,026 | 43,431,503 | 17,943,523 | 86,282,052 |

This is because expanding covered states/provinces in a country leads to an increase in population, which enlarges the problem (perhaps growing the deprivation costs). Thus, the more the supply frequency between distribution centers and vaccination centers, the less the growth rate of deprivation cost.

### 5.3. Sensitivity analysis

Sensitivity analysis examines how different values of an independent variable, known as input parameters, can affect the objective function and the obtained results. Also, the model is referred to as what-if analysis. Sensitivity analysis makes it easy to foresee the consequence of decisions given a specific range of variables and better understand it.

Sensitivity analysis is a common tool to analyze the goodness of the obtained solutions in a mathematical model; however, sensitivity analysis is impossible to perform for models having multiple uncertain parameters. But, in a robust optimization approach, we are making it possible to observe the sensitivity of the problem regarding its parameters and indices. In fact, this post-optimality analysis aims to discover the influence of input data fluctuations on the model and its solutions.

As previously mentioned, incorporating deprivation cost in the supply chain model is one of the most important contributions of this study. Figure 4 shows how deprivation cost varies by changing the size of the problem and periods. As already discussed, when the problem size increases, the growth rate of deprivation costs decreases. Actually, Figure 4 demonstrates that though deprivation costs are growing, the growth rate decreases.

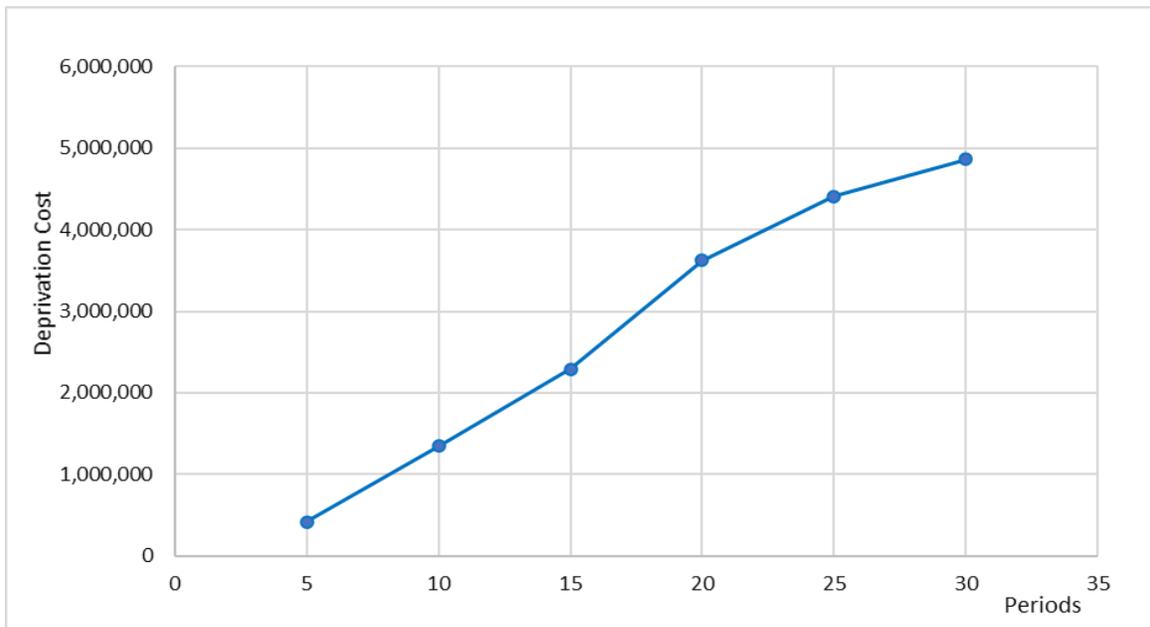

**Fig. 4.** Variation of deprivation cost for different problems.

In order to conduct sensitivity analysis, the two essential parameters of the proposed model, i.e., the maximum amount of allowed ordering for suppliers and the available budget for suppliers are considered. These experiments are designed based on 10% increase or decrease in nominal values of each parameter. However, their statuses are assessed when these values are unchanged. Obviously, if all of the considered parameters remain unchanged, it is redundant; accordingly, the status in which all of them are "unchanged" is removed, so 8 problem instances are considered.

**Table 4.** Results of sensitivity analysis on the value of budget and No. of orders.

| #Problem instance | Supplier Budget | Maximum of Orders | Objective function values | % Change in Objective function |
|---|---|---|---|---|
| 1 | Increased | Increased | 6,001,090 | 15.8% |
| 2 | Increased | Unchanged | 5,767,887 | 11.3% |
| 3 | Increased | Decreased | 5,576,142 | 7.6% |
| 4 | Unchanged | Increased | 5,488,043 | 5.9% |
| 5 | Unchanged | Decreased | 5,363,668 | 3.5% |
| 6 | Decreased | Increased | 5,322,210 | 2.7% |
| 7 | Decreased | Unchanged | 5,021,637 | -3.1% |
| 8 | Decreased | Decreased | 4,959,450 | -4.3% |

Table 4 shows that the percentage of change in the objective function complies with changes in the input parameters. For instance, in the first row of Table 4, increasing both of the supplier budget and the maximum order by 10% will cause a 15.8% increase in the objective function value. Instance #2 shows that if supplier budget increases, while the maximum of orders is unchanged, 11.3% increase concerning the linear version of instance #5 in Table 2 occurrs. Having compared the obtained solution of robust model and validated solutions of sensitivity analysis shows that robust optimization outperforms classic approaches and decreases computational time through changing parameters.

**Table 5.** Results of sensitivity analysis on the supplier budget.

| #Problem instance | Decrease in Suppliers Budget | Transportation Cost | Deprivation Cost | Objective function values |
|---|---|---|---|---|
| 1 | 5% | 2,020,438 | 2,213,987 | 5,210,537 |
| 2 | 10% | 1,845,775 | 2,264,201 | 5,021,637 |
| 3 | 15% | 1,738,433 | 2,298,679 | 4,898,344 |
| 4 | 20% | 1,588,274 | 2,322,617 | 4,712,805 |
| 5 | 25% | 1,344,697 | 2,391,234 | 4,471,770 |
| 6 | 30% | 1,116,991 | 2,438,653 | 4,227,910 |
| 7 | 40% | 936,789 | 2,518,739 | 4,056,213 |
| 8 | 50% | 523,183 | 2,635,148 | 3,661,404 |

Table 5 presents sensitivity analysis on suppliers' budgets, one of the critical parameters in the proposed model. As expected, when the available budget for suppliers decreases, it reduces transportation costs. It happens because of a reduction in supply frequency between distribution centers and vaccination centers. But surprisingly, a considerable increase in deprivation costs was observed. It could be interpreted in this way: as suppliers' budget decreases, the model tries to send as large as possible orders to supply vaccines. In other words, although the reduction in available budget causes a reduction in operating costs, such as transportation costs, the number of backorders and consequently deprivations costs would have a considerable increase, since deprivation costs relate to humanitarian issues and has more importance.